\title{Symmetries of quasiplatonic Riemann surfaces}

\author{Gareth A. Jones and David Singerman\\
School of Mathematics\\
University of Southampton\\
Southampton SO17  1BJ, U.K.\\
\\
Paul D.~Watson\\
Peter Symonds College\\
Winchester SO22 6RX, U.K.}

\documentclass[10pt]{article}
\usepackage{color}
\usepackage[latin1]{inputenc}
\usepackage{latexsym}
\usepackage{amsmath}
\usepackage{amssymb}
\usepackage{amsfonts}
\usepackage{tikz}
\newtheorem{thm}{Theorem}[section]
\newtheorem{lemma}[thm]{Lemma}

\date{}
\begin{document}

\maketitle

\begin{abstract}
We state and prove a corrected version of a theorem of Singerman, which relates the existence of symmetries (anticonformal involutions) of a quasiplatonic Riemann surface $\mathcal S$  (one uniformised by a normal subgroup $N$ of finite index in a cocompact triangle group $\Delta$) to the properties of the group $G=\Delta/N$. We give examples to illustrate the revised necessary and sufficient conditions for the existence of symmetries, and we relate them to properties of the associated dessins d'enfants, or hypermaps.
\end{abstract}

(secondary).

\medskip


\section{Introduction}

The category of compact Riemann surfaces is naturally equivalent to that of complex projective algebraic curves. Such a surface or curve $\mathcal S$ is real (defined over $\mathbb R$) if and only if it possesses an anticonformal involution called a {\sl symmetry}, in which case $\mathcal S$ is said to be {\sl symmetric}. (These specialised usages of these terms were introduced by Klein~\cite{Kle}; see~\cite{BCGG} for background.)

Among the most important Riemann surfaces $\mathcal S$ are the quasiplatonic surfaces, those uniformised by a normal subgroup $N$ of finite index in a cocompact triangle group 
\[\Delta=\Delta(l,m,n)=\langle X, Y, Z\mid X^l=Y^m=Z^n=XYZ=1\rangle;\]
this is equivalent to $\mathcal S$ carrying a regular dessin $\mathcal D$ in Grothendieck's terminology~\cite{Gro}, that is, an orientably regular hypermap in combinatorial language. (See~\cite{GG} or~\cite{JS} for the background to these connections.) In this case, $G:=\Delta/N$ can be identified with the automorphism group ${\rm Aut}\,{\mathcal D}$ of $\mathcal D$, the orientation-preserving automorphism group of the hypermap. In~\cite[Thm~2]{Sin74}, Singerman stated a result which relates the  existence of symmetries of such a surface $\mathcal S$ to the properties of $G$ (or equivalently of $\mathcal D$). As pointed out by Watson~\cite{Wat}, the stated conditions are sufficient for $\mathcal S$ to be symmetric, but are not necessary. A correct statement, based on~\cite[Thm~3.4]{Wat}, is as follows (the original~\cite{Sin74} omitted conditions~(3) and~(4)):

\begin{thm}\label{symsurf}
Let $\mathcal S$ be a quasiplatonic Riemann surface, uniformised by a normal subgroup $N$ of finite index in a cocompact triangle group $\Delta$, and let $x, y$ and $z$ be the images in $G:=\Delta/N$ of a canonical generating triple $X, Y, Z$ for $\Delta$. Then $\mathcal S$ is symmetric if and only if at least one of the following holds:
\begin{enumerate}
\item $G$ has an automorphism $\alpha: x\mapsto x^{-1}, \; y\mapsto y^{-1}$;
\item $G$ has an automorphism $\beta: x\mapsto y^{-1}, \; y\mapsto x^{-1}$ (possibly after a cyclic permutation of the canonical generators);
\item $\Delta$ has type $(2n,2n,n)$ for some $n$ (possibly after a cyclic permutation of its generators), $G$ has an automorphism $\gamma$ transposing $x$ and $y$, and the extension $\langle G,\gamma\rangle$ of $G$ by $\langle\gamma\rangle$ has an automorphism $\delta$ transposing $x$ and $x\gamma$; 
\item $\mathcal S$ has genus $1$.
\end{enumerate}
\end{thm}

The automorphisms $\alpha,\ldots,\delta$ in cases~(1), (2) and (3) must have order dividing $2$, and in cases~(2) and (3), $x$ and $y$ must have the same order. 

If we also assume that $\Delta$ is {\sl maximal\/} among all triangle groups normalising $N$ (these exist if $\mathcal S$ has genus $g>1$, but not if $g=1$), we obtain a version of this theorem which, although a little less general, is simpler to state and to prove, is closer to the original in~\cite{Sin74}, and for $g>1$ is equivalent to~\cite[Thm~1.5.10]{BCGG}, where $G$ is assumed to be the full group ${\rm Aut}^+{\mathcal S}$ of conformal automorphisms of $\mathcal S$. The authors are grateful to J\"urgen Wolfart for suggesting this alternative:

\begin{thm}\label{symsurf2}
Let $\mathcal S$, $N$, $\Delta$ and $G$ be as in Theorem~\ref{symsurf}, and suppose that $\Delta$ is maximal among all triangle groups containing $N$ as a normal subgroup. Then $\mathcal S$ is symmetric if and only if either
\begin{enumerate}
\item $G$ has an automorphism $\alpha: x\mapsto x^{-1}, \; y\mapsto y^{-1}$, or
\item $G$ has an automorphism $\beta: x\mapsto y^{-1}, \; y\mapsto x^{-1}$ (possibly after a cyclic permutation of the canonical generators).
\end{enumerate}
\end{thm}


\section{Combinatorial interpretation}\label{combint}

In view of the importance of quasiplatonic surfaces for the theories of dessins d'enfants and of maps and hypermaps, we will give combinatorial interpretations of the conditions in Theorem~\ref{symsurf}. Each quasiplatonic surface $\mathcal S$ inherits from $\Delta$ and $N$ a combinatorial structure  $\mathcal D$ called a {\sl regular dessin}, or {\sl orientably regular hypermap}, of type $(l, m, n)$. For our purposes, one can regard $\mathcal D$ as a triangulation of $\mathcal S$, the quotient by $N$ of the triangulation of the universal covering space $\hat{\mathcal S}$ naturally associated with $\Delta$, together with a preferred orientation of $\mathcal S$. The vertices of $\mathcal D$ can be coloured black, white or red, and termed {\sl hypervertices, hyperedges} or {\sl hyperfaces}, as they are quotients of fixed points of conjugates of $X$, $Y$ or $Z$, with valencies $2l, 2m$ or $2n$ respectively. Equivalently, in the language of dessins d'enfants, $\mathcal D$ is the inverse image, under a regular Bely\u\i\/ function $\beta:{\mathcal S}\to{\mathbb P}^1({\mathbb C})$, of the triangulation of the Riemann sphere ${\mathbb P}^1({\mathbb C})$ with vertices at the ramification points $0, 1$ and $\infty$ of $\beta$ coloured black, white and red, and with edges along $\mathbb R$.

The group $G=\Delta/N$ can be identified with the group ${\rm Aut}\,{\mathcal D}$ of automorphisms of $\mathcal D$; by definition this preserves the orientation and vertex colours of the triangulation, and since $\mathcal D$ is regular it acts regularly on incident pairs of vertices of given colours. We say that $\mathcal D$ is a {\sl reflexible\/} dessin (or a {\sl regular\/} hypermap) if the triangulation has an additional colour-preserving automorphism (induced by the extended triangle group $\Delta^*$ corresponding to $\Delta$) which reverses the orientation of $\mathcal S$, so that $\mathcal D$ is isomorphic to its mirror image $\overline{\mathcal D}$; otherwise $\mathcal D$ and $\overline{\mathcal D}$ form a {\sl chiral pair}. The {\sl black-white dual\/} ${\mathcal D}^{(01)}$ of $\mathcal D$ is the dessin formed by transposing the colours of the black and white vertices, or equivalently by replacing $\beta$ with $1-\beta$; similarly, there are black-red and white-red duals  ${\mathcal D}^{(02)}$ and ${\mathcal D}^{(12)}$ corresponding to $1/\beta$ and $\beta/(\beta-1)$.

For many purposes it is simpler, if less symmetric, to represent $\mathcal D$ by its {\sl Walsh map\/}~\cite{Wal} ${\mathcal W}=W({\mathcal D})$, a bipartite map on $\mathcal S$ formed by deleting the red vertices and their incident edges. This is the inverse image under $\beta$ of the unit interval $[0,1]\subset{\mathbb R}$. One can recover $\mathcal D$ as the stellation of $\mathcal W$, formed by placing a red vertex in each face of $\mathcal W$, and joining it by an edge to each incident black or white vertex. Many authors take this as the standard model of a dessin.

If one of the elliptic periods $l, m$ or $n$ of $\Delta$ is equal to $2$ (say $m=2$ without loss of generality), then the triangulation $\mathcal D$ is the barycentric subdivision $B({\mathcal M})$ of an orientably regular map $\mathcal M$ on $\mathcal S$, with vertices, edge-centres and face-centres at the black, white and red vertices of $\mathcal D$, and with orientation-preserving automorphism group ${\rm Aut}^+{\mathcal M}=G$. In the notation of Coxeter and Moser~\cite[Ch.~8]{CM}, $\mathcal M$ has type $\{n, l\}$, meaning that vertices and faces have valencies $l$ and $n$. Conversely, any orientably regular map $\mathcal M$ on $\mathcal S$ determines a dessin ${\mathcal D}=B({\mathcal M})$ with $m=2$. The black-red dual ${\mathcal D}^{(02)}$ of $\mathcal D$ corresponds to the vertex-face dual map ${\mathcal M}^*$ of $\mathcal M$.

In Theorem~\ref{symsurf}, condition~(1) corresponds to $\mathcal D$ being reflexible, as happens whenever $\mathcal S$ has genus $0$, for example. Condition~(2) corresponds to $\mathcal D$ being isomorphic to $\overline{{\mathcal D}^{(01)}}$, the mirror image of its black-white dual, as happens (after a cyclic permutation of generators) for the Edmonds maps ${\mathcal M}\cong\overline{{\mathcal M}^*}$ of genus $7$ and type $\{7,7\}$; these are a chiral pair of orientably regular embeddings of the complete graph $K_8$ in Macbeath's curve~\cite{Macb65}, described by Coxeter in~\cite[\S21.3]{Cox} and denoted by C7.2 in Conder's list of chiral maps~\cite{Con} (see also~\cite[p.~29]{Sin74}. Neither of these conditions applies to the chiral dessins of genus $1$, such as the embeddings ${\mathcal M}=\{4,4\}_{1,2}$ and $\overline{\mathcal M}=\{4,4\}_{2,1}$ of $K_5$, or the embeddings $\{3,6\}_{1,2}$ and $\{3,6\}_{2,1}$ of $K_7$, with $G\cong AGL_1(5)$ or $AGL_1(7)$ (see~\cite[Ch.~8]{CM}): these groups have only inner automorphisms, and none satisfying~(1) or (2), but the underlying tori, uniformised by square and hexagonal lattices respectively, admit symmetries, so condition~(4) is required in Theorem~\ref{symsurf}. 

Condition~(3) is a little harder to explain. The group $G={\rm Aut}\,{\mathcal D}$ acts regularly on edges of the Walsh map $\mathcal W$ of $\mathcal D$; the existence of an automorphism $\gamma$ as in~(3) is equivalent to $\mathcal W$, regarded as an uncoloured map, being orientably regular, with a half-turn $\gamma$ reversing an edge, so that ${\rm Aut}^+{\mathcal W}=\langle G,\gamma\rangle$ acts regularly on arcs (directed edges) of $\mathcal W$; the existence of $\delta$ is equivalent to ${\mathcal W}\cong\overline{{\mathcal W}^*}$. In \S\ref{exceptionexample} we will give examples of dessins $\mathcal D$ satisfying~(3) but not (1), (2) or (4), with underlying Riemann surfaces $\mathcal S$ admitting symmetries, so that this condition is required for a correct statement of Theorem~\ref{symsurf}. (Note that $G$ is not normal in $\langle G,\gamma, \delta\rangle$ since $x$ is conjugate to $\gamma x\not\in G$.)

Dessins $\mathcal D$ satisfying~(3) or (4), but not (1) or (2), and hence not covered by~\cite[Theorem~2]{Sin74}, all correspond to triangle groups with periods $2n, 2n, n>2$, or have genus $1$. Thus the original statement is correct when restricted to {\sl maps\/} of genus $g\ne 1$, and indeed when applied to dessins of all types except permutations of $(2,3,6), (3,3,3)$ and $(2n,2n,n)$ for $n\ge 2$.

\section{Groups containing triangle groups}

In order to prove Theorem~\ref{symsurf} we need to consider which isometry groups $\tilde\Delta$ (of $\hat{\mathcal S}={\mathbb P}^1({\mathbb C})$, $\mathbb C$ or $\mathbb H$) can contain a cocompact triangle group $\Delta=\Delta(l,m,n)$ as a subgroup of index $2$. Here are some examples:

\begin{enumerate}
\item The most obvious example is the extended triangle group $\Delta^*=\Delta^*(l,m,n)$, the extension of $\Delta$ by a reflection $T$ in the geodesic through the fixed points of two of its canonical generators $X, Y, Z$, which it inverts by conjugation; each choice of a pair from the generating triple gives the same group $\Delta^*$.

\item If two of $l, m$ and $n$ are equal, say $l=m$ without loss of generality, $\Delta$ has index $2$ in a triangle group $\Delta^==\Delta(l,2,2n)$, the extension of $\Delta$ by a rotation of order $2$ transposing the fixed points of $X$ and $Y$, and transposing $X$ and $Y$ by conjugation; $\Delta$ also has index $2$ in a group $\Delta^{\times}=\Delta^{\times}(l,m,n)$, the extension of $\Delta$ by a reflection transposing the fixed points of $X$ and $Y$, and sending $X$ to $Y^{-1}$ and $Y$ to $X^{-1}$ by conjugation.

\item If $l=m=n$ then in addition to $\Delta^*$ we obtain three triangle groups $\Delta^=$ and three groups $\Delta^{\times}$ as in case~(2), depending in the choice of a pair from the canonical generating triple $X, Y, Z$.

\end{enumerate}

In the standard notation for NEC groups introduced by Macbeath~\cite{Macb67}, $\Delta$ has signature $(0, +, [l, m, n], \{\,\})$, while $\Delta^*$, $\Delta^=$ and $\Delta^{\times}$ have signatures $(0, +, [\;], \{(l, m, n)\})$,  $(0, +, [l, 2, 2n], \{\,\})$ and  $(0, +, [l], \{(n)\})$ respectively. The following result shows that these groups are in fact the only possibilities for $\tilde\Delta$:

\begin{lemma}\label{index2}
Let $\Delta$ be a cocompact triangle group $\Delta(l,m,n)$, and let $\tilde\Delta$ be an isometry group containing $\Delta$ as a subgroup of index $2$. Then either 
\begin{enumerate}
\item $l, m$ and $n$ are distinct, and $\tilde\Delta=\Delta^*$, or
\item just two of $l, m$ and $n$ are equal, and $\tilde\Delta=\Delta^*$, $\Delta^=$ or $\Delta^{\times}$, or
\item $l=m=n$, and $\tilde\Delta=\Delta^*$, one of the three groups $\Delta^=$, or one of the three groups $\Delta^{\times}$.
\end{enumerate}
\end{lemma}

\noindent{\sl Proof.} Any index $2$ inclusion must be normal, so the groups $\tilde\Delta$ containing $\Delta$ as a subgroup of index $2$ are contained in the normaliser $N(\Delta)$ of $\Delta$ in the full isometry group, and correspond to the involutions in $N(\Delta)/\Delta$. In case~(1) we have $N(\Delta)=\Delta^*$, with $N(\Delta)/\Delta\cong C_2$, so $\tilde\Delta=\Delta^*$. In case~(2) we can assume without loss of generality that $l=m\ne n$, so $N(\Delta)=\Delta^*(l,2,2n)$ with $N(\Delta)/\Delta\cong V_4$; the three involutions in this group yield the possibilities $\tilde\Delta=\Delta^*$, $\Delta^=$ and $\Delta^{\times}$. In case~(3) we have $N(\Delta)=\Delta^*(2,3,2n)$ and $N(\Delta)/\Delta\cong \Delta^*(2,3,2)\cong S_3\times C_2$; there are seven involutions in this group (the central involution, generating the direct factor $C_2$, and two conjugacy classes of three non-central involutions), giving the seven subgroups $\tilde\Delta$ listed in~(3). \hfill$\square$

\medskip

Now suppose that $G=\Delta/N$ for some normal subgroup $N$ of $\Delta$. If $N$ is also normal in one of the groups $\tilde\Delta=\Delta^*$, $\Delta^=$ or $\Delta^{\times}$ containing $\Delta$ with index $2$, let $\tilde G=G^*$, $G^=$ or $G^{\times}$ be the corresponding extension $\tilde\Delta/N$ of $G$ by an involution acting as above on the canonical generators $x, y, z$ of $G$. Now $N$ corresponds to a regular dessin $\mathcal D$, with $G\cong{\rm Aut}\,{\mathcal D}$, and normality of $N$ in $\Delta^*$ corresponds to $\mathcal D$ being reflexible, that is, ${\mathcal D}\cong\overline{\mathcal D}$, with $G^*$ the full automorphism group of the hypermap; normality in $\Delta^=$ corresponds to ${\mathcal D}\cong{\mathcal D}^{(01)}$, with $G^=={\rm Aut}^+{\mathcal W}$, where ${\mathcal W}=W({\mathcal D})$ is regarded as an orientably regular uncoloured map; finally, normality of $N$ in $\Delta^{\times}$ corresponds to $\mathcal D$ being isomorphic to the mirror image of its dual ${\mathcal D}^{(01)}$. (When $l=m=n$ there are three duals to consider.)

\section{Proof of Theorems~\ref{symsurf} and~\ref{symsurf2}}

We will now prove Theorem~\ref{symsurf}, with a short digression in \S4.2 to deal with Theorem~\ref{symsurf2}.

\subsection{The conditions are sufficient}

Let $\mathcal D$ have type $(l, m, n)$, so it corresponds to a normal subgroup $N$ of $\Delta=\Delta(l,m,n)$, with $\Delta/N\cong G={\rm Aut}\,{\mathcal D}$. We will show that each of conditions~(1) to (4) in Theorem~\ref{symsurf} is each sufficient for $\mathcal S$ to admit a symmetry.

If condition~(1) holds, let $G^*$ be the extension of $G$ by $\langle\alpha\rangle\cong C_2$, with $\alpha$ acting naturally by conjugation on $G$. The epimorphism $\Delta\to G$ with kernel $N$ extends to an epimorphism $\Delta^*=\langle \Delta,T\rangle\to G^*,\; T\mapsto\alpha$ with kernel $N$, where $\Delta^*=\Delta^*(l, m, n)$, and the reflection $T$ induces a symmetry of $\mathcal S$. A similar argument applies to condition~(2), with $\Delta^{\times}$ used instead of $\Delta^*$.

If condition~(3) holds, the epimorphism $\Delta\to G$ extends, firstly to an epimorphism $\Delta^=\to\langle G,\gamma\rangle=G^=$, and then to an epimorphism $(\Delta^=)^{\times}\to \langle G, \gamma, \delta\rangle=(G^=)^{\times}$, with $\delta$ lifting to a reflection inducing a symmetry of $\mathcal S$.

If condition~(4) holds, $\mathcal S$ is a torus ${\mathbb C}/\Lambda$ for some lattice $\Lambda\subset{\mathbb C}$. Having genus $1$, $\mathcal D$ must have a period $k=l, m$ or $n$ greater than $2$, so $\mathcal S$ admits an automorphism of order $k>2$ with a fixed point. This lifts to a rotation of order $k>2$ of $\Lambda$, so $\Lambda$ is a square or hexagonal lattice. In either case $\Lambda$ admits reflections of $\mathbb C$, inducing symmetries of $\mathcal S$.

\subsection{The conditions are necessary}

For the converse, suppose that $\mathcal S$ is symmetric. If $\mathcal S$ has genus $g=0$ then condition~(1) holds, and condition~(4) deals with the case $g=1$, so we may assume that $g\ge 2$. Any symmetry of $\mathcal S$ is induced by an orientation-reversing isometry $R$ of $\mathbb H$ normalising $N$ and satisfying $R^2\in N$. Let $\Gamma=\langle\Delta, R\rangle$. Since $\Delta$ and $R$ normalise $N$, $\Gamma$ is contained in the normaliser of $N$ in the isometry group $PGL_2({\mathbb R})$. This is an NEC group, and hence so is $\Gamma$. Since it contains $R$, $\Gamma$ is a proper NEC group, so let $\Gamma^+$ be its orientation-preserving Fuchsian subgroup of index $2$.

If $\Delta$ is a normal subgroup of $\Gamma$ then $|\Gamma:\Delta|=2$, so $\Delta=\Gamma^+$. Then Lemma~\ref{index2} shows that $\Gamma$ is either $\Delta^*$, giving condition~(1), or a group $\Delta^{\times}$, giving condition~(2). We may therefore assume that we have a non-normal inclusion $\Delta<\Gamma$. Since $\Gamma^+$ is a Fuchsian group containing the triangle group $\Delta$, it is also a triangle group, as shown by Singerman in~\cite{Sin72}; moreover, $\Delta$ is a proper subgroup of $\Gamma^+$ since it is not normal in $\Gamma$. This immediately proves Theorem~1.2, since it contradicts the maximality of $\Delta$ assumed there. Continuing with the proof of Theorem~\ref{symsurf}, the inclusion  $\Delta<\Gamma^+$ must appear in Singerman's list of triangle group inclusions~\cite{Sin72}, shown in Table~1. There are seven sporadic examples and seven infinite families: cases (a), (b) and (c) are normal inclusions, while cases (A) to (K) are non-normal. In the fifth column, $P$ denotes the permutation group induced by $\Gamma^+$ on the cosets of $\Delta$.

\begin{table}[ht]
\centering
\begin{tabular}{| p{0.8cm} | p{1.6cm} | p{1.8cm} | p{1.3cm} | p{2.3cm} | p{1.2cm} |}
\hline
Case & Type of $\Delta$ & Type of $\Gamma^+$ & $|\Gamma^+:\Delta|$ & $P$ & Thm~1.1 \\
\hline\hline
a & $(s,s,t)$ & $(2,s,2t)$ & $2$ & $S_2$ & (1), (2) \\
\hline
b & $(t,t,t)$ & $(3,3,t)$ & $3$ & $A_3$ & (2) \\
\hline
c & $(t,t,t)$ & $(2,3,2t)$ & $6$ & $S_3$ & (1), (2) \\
\hline\hline
A & $(7,7,7)$ & (2,3,7) & $24$ & $L_2(7)$ & (1) \\
\hline
B & $(2,7,7)$ & (2,3,7) & $9$ & $L_2(8)$ & (2) \\
\hline
C & $(3,3,7)$ & (2,3,7) & $8$ & $L_2(7)$ & (2) \\
\hline
D & $(4,8,8)$ & (2,3,8) & $12$ & $(C_4\times C_4)\rtimes S_3$ & (2) \\ 
\hline
E & $(3,8,8)$ & (2,3,8) & $10$ & $PGL_2(9)$ & (2) \\
\hline
F & $(9,9,9)$ & (2,3,9) & $12$ & $L_2({\mathbb Z}_9)$ & (2) \\
\hline
G & $(4,4,5)$ & (2,4,5) & $6$ & $S_5$ & (2) \\
\hline
H & $(n,4n,4n)$ & $(2,3,4n)$ & $6$ & $S_4$ & (2) \\
\hline
I & $(n,2n,2n)$ & $(2,4,2n)$ & $4$ & $D_4$ & (2) \\
\hline
J & $(3,n,3n)$ & $(2,3,2n)$ & $4$ & $A_4$ & (1) \\
\hline
K & $(2,n,2n)$ & $(2,3,2n)$ & $3$ &  $S_3$ & (1) \\
\hline

\end{tabular}
\caption{Inclusions between Fuchsian triangle groups.}
\label{table:MobSuz}
\end{table}

In all cases except (a), with $s=2t$, and (b), $\Gamma^+$ has no repeated periods, so Lemma~\ref{index2} shows that $\Gamma$ is the extended triangle group $(\Gamma^+)^*$. As shown by Watson in~\cite[Appendix]{Wat}, in each such case the dessin corresponding to the inclusion $\Delta<\Gamma^+$ is reflexible (see \S4.3 for an example), so $\Delta$ has index $2$ in a proper NEC group $\tilde\Delta\le\Gamma$. Lemma~\ref{index2} shows that $\tilde\Delta=\Delta^*$ or $\Delta^{\times}$. Since $\tilde\Delta$, as a subgroup of $\Gamma$, normalises $N$, it follows that $\tilde G:=\tilde\Delta/N$ is an extension of $G=\Delta/N$ by an automorphism satisfying condition~(1) or (2) of Theorem~\ref{symsurf}  respectively, as indicated in the final column of Table~1.

In the two exceptional cases, a pair of repeated periods allows the additional possibility that $\Gamma=(\Gamma^+)^{\times}$ (not $(\Gamma^+)^=$, since $\Gamma$ is a proper NEC group). In case (b), if $\Gamma=(\Gamma^+)^{\times}$ then $\Delta$ is normal in $\Gamma$ with quotient $C_6$, and the involution in $\Gamma/\Delta$ corresponds to a proper NEC group $\tilde\Delta<\Gamma$ containing $\Delta$ with index $2$, as before. Hence there remains only case (a) with $s=2t$, or equivalently $\Delta=\Delta(2n,2n,n)$, $\Gamma^+=\Delta(2n,2,2n)$ and $\Gamma=\Delta^{\times}(2n,2,2n)$, considered in \S\ref{exception}.

\subsection{Example of a typical case}

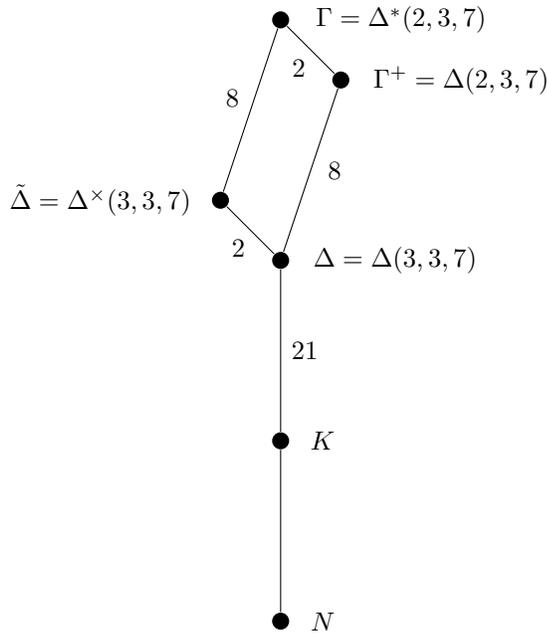
\begin{figure}[h!]
\label{237lattice}
\begin{center}
 \begin{tikzpicture}[scale=0.8, inner sep=0.8mm]

\node (A) at (0,10) [shape=circle, fill=black] {};
\node (B) at (1,9) [shape=circle, fill=black] {};
\node (C) at (-1,7) [shape=circle, fill=black] {};
\node (D) at (0,6) [shape=circle, fill=black] {};
\node (E) at (0,3) [shape=circle, fill=black] {};
\node (F) at (0,0) [shape=circle, fill=black] {};

\draw (A) to (B);
\draw (A) to (C);
\draw (B) to (D);
\draw (C) to (D);
\draw (D) to (E);
\draw (E) to (F);

\node at (2,10) {$\Gamma=\Delta^*(2,3,7)$};
\node at (3,9) {$\Gamma^+=\Delta(2,3,7)$};
\node at (-3,7) {$\tilde\Delta=\Delta^{\times}(3,3,7)$};
\node at (1.9,6) {$\Delta=\Delta(3,3,7)$};
\node at (0.7,3) {$K$};
\node at (0.7,0) {$N$};

\node at (0.3,9.2) {$2$};
\node at (-0.7,6.2) {$2$};
\node at (-0.8,8.7) {$8$};
\node at (0.9,7.5) {$8$};
\node at (0.4,4.5) {$21$};

\end{tikzpicture}

\end{center}
 \caption{Lattice of some subgroups of $\Delta^*(2,3,7)$} 
\end{figure}

Suppose that the inclusion of triangle groups $\Delta<\Gamma^+$ is case~C in Table~1. Thus $\Delta=\Delta(3,3,7)$ is a subgroup of index $8$ in $\Gamma^+=\Delta(2,3,7)$, so that $\Gamma=\Delta^*(2,3,7)$ by Lemma~\ref{index2}. Since $N$ is a subgroup of $\Delta$ and is normal in $\Gamma$, it is contained in the core $K$ of $\Delta$ in $\Gamma$. This is a normal subgroup of $\Gamma$, a surface group of genus $3$ contained in $\Gamma^+$, with $\Gamma^+/K\cong P=L_2(7)$ and $\Gamma/K\cong PGL_2(7)$; the Riemann surface ${\mathbb H}/K$ uniformised by $K$ is Klein's quartic curve $x^3y+y^3z+z^3x=0$. The image of $\Delta$ in $PGL_2(7)$ is a subgroup $H$ of order $21$ and index $8$ in $L_2(7)$: this is the stabiliser of a point in the natural action on the projective line ${\mathbb P}^1({\mathbb F}_7)$, isomorphic to the unique subgroup of index $2$ in $AGL_1(7)$. The stabiliser in $PGL_2(7)$ of this point is isomorphic to $AGL_1(7)$, an extension of $H$ by an involution which acts on its two canonical generators of order $3$ as in condition~(2) of Theorem~\ref{symsurf}. This lifts to a proper NEC subgroup $\tilde\Delta$ of $\Gamma$ which contains $\Delta$ with index $2$, acting in the same way on its two generators of order $3$. By Lemma~\ref{index2} $\tilde\Delta$ must be $\Delta^{\times}=\Delta^{\times}(3,3,7)$, so it contains a reflection which, since it normalises $N$, induces on $G:=\Delta/N$ an automorphism satisfying condition~(2). Figure~1 shows the inclusions between these subgroups of $\Gamma$; edges are labelled with indices of inclusions. The dessin corresponding to the inclusion $\Delta<\Gamma^+$ is shown in Figure~2 as a map on the sphere, where we have changed generators to take $\Gamma^+=\Delta(3,2,7)$; the generators $x, y$ and $z$ of $G$ of order $3, 3$ and $7$ correspond, as in~\cite[Thm~1]{Sin70}, to short cycles of the elliptic generators of $\Gamma^+$, at the two vertices and the one face of valency $1$; the obvious reflection transposes and inverts $x$ and $y$.

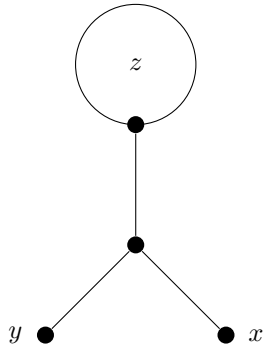
\begin{figure}[h!]
\label{237map}
\begin{center}
 \begin{tikzpicture}[scale=0.8, inner sep=0.8mm]

\node (A) at (0,3.5) [shape=circle, fill=black] {};
\node (B) at (0,1.5) [shape=circle, fill=black] {};
\node (C) at (-1.5,0) [shape=circle, fill=black] {};
\node (D) at (1.5,0) [shape=circle, fill=black] {};

\draw (A) to (B);
\draw (B) to (C);
\draw (B) to (D);

\draw (1,4.5) arc (0:360:1);

\node at (2,0) {$x$};
\node at (-2,0) {$y$};
\node at (0,4.5) {$z$};

\end{tikzpicture}

\end{center}
 \caption{Map for the inclusion $\Delta(3,3,7)<\Delta(2,3,7)$} 
\end{figure}

For specific instances of such subgroups $N$ one can use Macbeath's construction in~\cite{Macb61} of an infinite sequence of Hurwitz groups $\Gamma^+/N$: for any integer $r\ge 1$, the group $N=K'K^r$ generated by the commutators and $r$th powers of elements of $K$ is a characteristic subgroup of $K$ and hence a normal subgroup of $\Gamma$, of index $336r^6$. The surface ${\mathcal S}={\mathbb H}/N$ carries a regular chiral dessin $\mathcal D$ of type $(3,3,7)$ and genus $1+2r^6$, a regular covering of the chiral dessin of type $(3,3,7)$ and genus $1$ corresponding to the inclusion $K<\Delta$. The automorphism group of the Riemann surface $\mathcal S$ is an extension $\Gamma^+/N$ of an abelian group $K/N\cong C_r^6$ by $\Gamma^+/K\cong L_2(7)$; if we include anti-conformal automorphisms we obtain an extension $\Gamma/N$ of $K/N$ by $\Gamma/K\cong PGL_2(7)$.

\subsection{The exceptional case}\label{exception}

\begin{figure}[h!]
\label{lattice}
\begin{center}
 \begin{tikzpicture}[scale=0.8, inner sep=0.8mm]

\node (A) at (0,8) [shape=circle, fill=black] {};
\node (B) at (-2,6) [shape=circle, fill=black] {};
\node (C) at (0,6) [shape=circle, fill=black, label=right: $\Gamma$] {};
\node (D) at (2,6) [shape=circle, fill=black] {};
\node (E) at (-4,4) [shape=circle, draw] {};
\node (F) at (-2,4) [shape=circle, draw, label=right: $\Delta^{\times}$] {};
\node (G) at (0,4) [shape=circle, fill=black] {};
\node (H) at (-2,2) [shape=circle, draw] {};
\node (I) at (2,2) [shape=circle, draw] {};
\node (J) at (0,0) [shape=circle, fill=black,] {};

\draw (A) to (E);
\draw (A) to (G);
\draw (A) to (D);
\draw (B) to (F);
\draw (B) to (I);
\draw (D) to (H);
\draw (E) to (H);
\draw (F) to (H);
\draw (H) to (J);
\draw (I) to (J);

\node at (2,8) {$\Lambda=\Delta^*(4,2,2n)$};
\node at (-0.8,6) {$(\Delta^=)^{\times}$};
\node at (-4.5,6) {$N(\Delta)=\Delta^*(2n,2,2n)$};
\node at (4,6) {$\Lambda^+=\Delta(4,2,2n)$};
\node at (-6.2,4) {$\Delta^*=\Delta^*(2n,2n,n)$};
\node at (3,4) {$\Gamma^+=\Delta^==\Delta(2n,2,2n)$};
\node at (-4,2) {$\Delta=\Delta(2n,2n,n)$};
\node at (4.2,2) {$\Delta^L\;(L\in\Lambda\setminus N(\Delta))$};
\node at (1.7,0) {$K=\Delta\cap\Delta^L$};

\end{tikzpicture}

\end{center}
 \caption{Lattice of groups in the exceptional case} 
\end{figure}
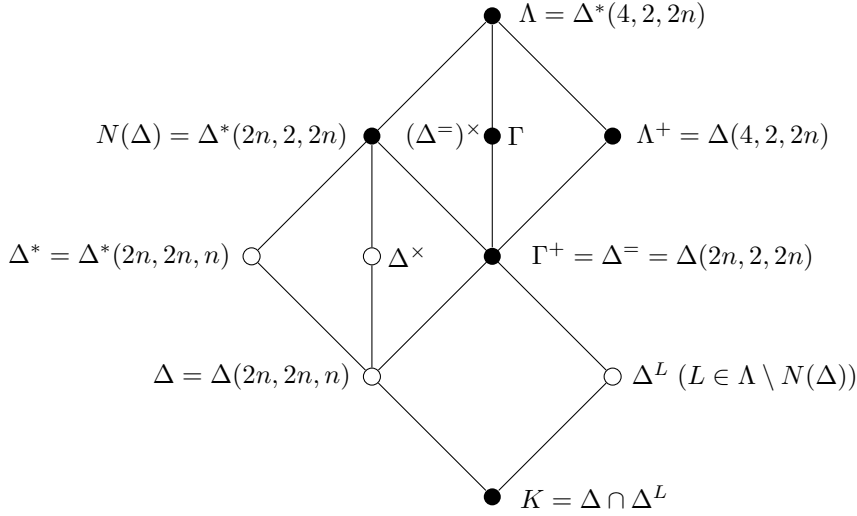

The exceptional case in the proof of Theorem~\ref{symsurf} arises in case~(a) of Table~1 when $s=2t$. Putting $t=n$ we have $\Delta=\Delta(2n,2n,n)$, $\Gamma^+=\Delta^==\Delta(2n,2,2n)$ and $\Gamma=\Delta^{\times}(2n,2,2n)$. In this case there is no proper NEC group $\tilde\Delta\le\Gamma$ containing $\Delta$ with index $2$: indeed, the only proper NEC groups containing $\Delta$ with index $2$ are $\Delta^*$ and $\Delta^{\times}$, both of which are subgroups of $\Delta^*(2n,2,2n)$ rather than of $\Delta^{\times}(2n,2,2n)$.

The inclusions between these NEC groups are shown in Figure~3, where $\Lambda$ is the maximal NEC group $\Delta^*(4,2,2n)$ containing $\Gamma$. Black and white vertices indicate normal and non-normal subgroups of $\Lambda$, and edges indicate inclusions, all of index $2$. The normaliser $N(\Delta)$ of $\Delta$ (in the isometry group of $\mathbb H$) is the extended triangle group $\Delta^*(2n,2,2n)$, which has index $2$ in $\Lambda$. The core of $\Delta$ in $\Lambda$ is the group $K=\Delta\cap\Delta^L$, where $L$ is any element of $\Lambda\setminus N(\Delta)$; this is a quadrilateral group $\Delta(n,n,n,n)$. The quotient $\Lambda/K$ is isomorphic to $\Delta^*(4,2,2)\cong D_4\times C_2$, the automorphism group of the regular map $\{2,4\}$ on the sphere; this can be seen by applying the natural epimorphism $\Lambda=\Delta^*(4,2,2n)\to\Delta^*(4,2,2)$, with kernel $K$.

To complete the proof, since $N$ is normal in $\Gamma$ with $G\cong\Delta/N$, and since $\Delta\le\Delta^=\le\Gamma$, there is an extension $G^==\langle G,\gamma\rangle\cong\Delta^=/N$ of $G$, where $\gamma$ is an automorphism of order $2$ of $G$ transposing its canonical generators $x$ and $y$ of order $2n$. The canonical generators of $G^=$, the images of those of $\Delta^==\Delta(2n,2,2n)$, are $x, \gamma$ and $\gamma x^{-1}$. Since $\Delta^=\le\Gamma=(\Delta^=)^{\times}$, there is an extension $(G^=)^{\times}=\langle G^=,\delta\rangle\cong\Gamma/N$ of $G^=$, where $\delta$ is an automorphism of order $2$ of $G^=$ sending $x$ to $(\gamma x^{-1})^{-1}=x\gamma$. Thus condition~(3) is satisfied. \hfill$\square$

\subsection{Example of the exceptional case}\label{exceptionexample}

In the exceptional case, $N$ is a torsion-free subgroup of finite index in $\Delta=\Delta(2n,2n,n)$; it is normal $\Gamma=(\Delta^=)^{\times}=\Delta^{\times}(2n,2,2n)$, and therefore contained in $K$. Any torsion-free characteristic subgroup $N$ of finite index in $K$ will correspond to a dessin $\mathcal D$ with $G={\rm Aut}\,{\mathcal D}$ satisfying condition~(3) of Theorem~\ref{symsurf}: the normality of $N$ in $\Gamma^+$ and $\Gamma$ provides the required automorphisms $\gamma$ and $\delta$. However, such a subgroup $N$ is normal in $\Lambda$, and hence in $\Delta^*$ and $\Delta^{\times}$, so~(1) and (2) are also satisfied. To show that condition (3) is independent of the others, and hence needed for a correct statement of Theorem~\ref{symsurf}, we will give an example where $N$ is normal in $\Gamma$ but not in $\Lambda$, with $n>2$ so that $\mathcal S$ does not satisfy~(4). If conditions~(1) or (2) were satisfied then $N$ would be normal in $\Delta^*$ or $\Delta^{\times}$ as well as in $\Gamma$, and hence normal in $\Lambda$, which is false.

In~\cite{Big}, Biggs showed that the complete graph $K_q$ on $q$ vertices has an orientably regular embedding if and only if $q$ is a prime power. When $q=2^e$ the examples he constructed are maps ${\mathcal M}_1$ of type $\{q-1,q-1\}$ and genus $(q-1)(q-4)/4$ with ${\rm Aut}^+{\mathcal M}_1\cong AGL_1(q)$; for instance, when $q=8$ they are the Edmonds maps. Each such map corresponds to a normal subgroup $M_1$ of $\Delta(n,2,n)$ with $\Delta(n,2,n)/M_1\cong AGL_1(q)$, where $n=q-1$. Using the fact that ${\rm Aut}\,AGL_1(q)=A\Gamma L_1(q)$, James and Jones~\cite{JJ} showed that if $q=2^e\ge 8$ then ${\mathcal M}_1$ is chiral and $\overline{\mathcal M}_1\cong{\mathcal M}_1^*$. 

We need these two properties to be satisfied by a bipartite map, which can then be the Walsh map of a dessin $\mathcal D$ as in the combinatorial explanation of condition~(3) in \S\ref{combint}. Since ${\mathcal M}_1$ is not bipartite, we construct a bipartite covering of it with the same properties. Let ${\mathcal M}_2$ be the orientably regular map of type $\{2,2\}$ on the sphere, with two vertices, joined by two edges, so that ${\rm Aut}^+{\mathcal M}_2\cong V_4$. The join ${\mathcal M}_1\vee{\mathcal M}_2$ of these two maps is an orientably regular map ${\mathcal M}_3$ of type $\{2n,2n\}$, corresponding to a torsion-free normal subgroup $N=N_1\cap N_2$ of finite index in $\Gamma^+=\Delta(2n,2,2n)$, where $N_1$ and $N_2$ are the inverse images of $M_1$ and $1$ in $\Gamma^+$ under the natural epimorphisms $\Delta(2n,2,2n)\to\Delta(n,2,n)$ and $\Delta(2n,2,2n)\to\Delta(2,2,2)\cong V_4$. Since the groups $\Gamma^+/N_i\cong AGL_1(q)$ and $V_4$ for $i=1, 2$ have no non-trivial common quotients, we have $\Gamma^+=N_1N_2$; thus
${\mathcal M}_3={\mathcal M}_1\times{\mathcal M}_2$ and
\[{\rm Aut}^+{\mathcal M}_3\cong\Gamma^+/N\cong (N_2/N)\times(N_1/N)\cong AGL_1(q)\times V_4.\]
These subgroups are shown in Figure~4, with black and white vertices indicating normal and non-normal subgroups of $\Lambda$.

Now $K$ is the unique normal subgroup of $\Gamma^+$ with quotient group $V_4$, so $N_2=K$ and hence $N\le\Delta$. Since ${\mathcal M}_i\cong\overline{{\mathcal M}_i^*}$ for $i=1,2$, each $N_i$ is normal in $\Gamma$, and hence so is $N$, so ${\mathcal M}_3\cong\overline{{\mathcal M}_3^*}$. If $N$ were normal in $\Lambda$ then it would be normal in $\Delta(2n,2,2n)^*$; since the direct factor $N_1/N\cong V_4$ is a characteristic subgroup of $\Gamma^+/N$ (as the centraliser of the subgroup $N_2/N\cong AGL_1(q)$ generated by the elements of odd order) it would follow that $N_1$ is normal in $\Delta(2n,2,2n)^*$, contradicting the chirality of ${\mathcal M}_1$. Thus $N$ is not normal in $\Lambda$, as required.

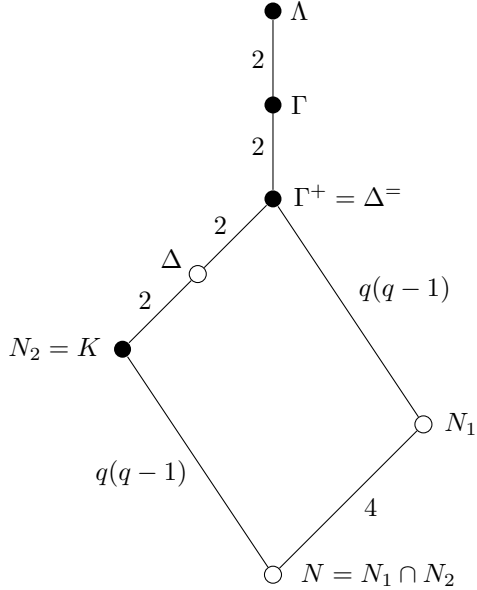
\begin{figure}[h!]
\label{examplelattice}
\begin{center}
 \begin{tikzpicture}[scale=0.5, inner sep=0.8mm]

\node (A) at (0,9) [shape=circle, fill=black] {};
\node (C) at (0,6.5) [shape=circle, fill=black] {};
\node (G) at (0,4) [shape=circle, fill=black] {};
\node (H) at (-2,2) [shape=circle, draw] {};
\node (J) at (-4,0) [shape=circle, fill=black,] {};
\node (L) at (4,-2) [shape=circle, draw] {};
\node (M) at (0,-6) [shape=circle, draw] {};

\draw (A) to (G);
\draw (G) to (L);
\draw (G) to (H);
\draw (H) to (J);
\draw (J) to (M);
\draw (L) to (M);

\node at (0.7,9) {$\Lambda$};
\node at (0.7,6.5) {$\Gamma$};
\node at (2,4.1) {$\Gamma^+=\Delta^=$};
\node at (-2.7,2.4) {$\Delta$};
\node at (-5.8,0) {$N_2=K$};
\node at (5,-2) {$N_1$};
\node at (2.8,-6) {$N=N_1\cap N_2$};

\node at (-0.4,7.7) {$2$};
\node at (-0.4,5.4) {$2$};
\node at (-1.4,3.3) {$2$};
\node at (-3.4,1.3) {$2$};
\node at (3.5,1.6) {$q(q-1)$};
\node at (-3.5,-3.3) {$q(q-1)$};
\node at (2.6,-4.2) {$4$};

\end{tikzpicture}

\end{center}
\caption{Lattice of groups in the example of the exceptional case} 
\end{figure}

This construction can be interpreted combinatorially as follows. As a covering of the bipartite map ${\mathcal M}_2$, the map ${\mathcal M}_3$ is also bipartite, so it is the Walsh map ${\mathcal W}=W({\mathcal D})$ of a dessin $\mathcal D$: this is a regular dessin of type $(2n,2n,n)$, corresponding to the normal inclusion of $N$ in $\Delta=\Delta(2n,2n,n)$. Condition~(3) corresponds to the fact that $\mathcal W$ is orientably regular and ${\mathcal W}\cong\overline{{\mathcal W}^*}$. The failure of conditions~(1) and (2) corresponds to $\mathcal D$ not being isomorphic to $\overline{\mathcal D}$ or $\overline{{\mathcal D}^{(01)}}$.

In this example, $\mathcal S$ has genus $g=n^2-n-1$. This is minimised when $e=3$, so $g=41$ and the chiral maps $\mathcal W$ and ${\mathcal W}^*$ correspond to C41.24 in~\cite{Con}. 

Further examples of this type can be found by using the `Macbeath trick'~\cite{Macb61} as in \S4.3. Let $r$ be any integer coprime to $2n$. Since $N'N^r$ is a characteristic subgroup of $N$, it is normal in $\Gamma$. If $N'N^r$ were normal in $\Lambda$ then, since $N/N'N^r$ is a normal subgroup of $\Lambda/N'N^r$ (being generated by its elements of order $r$), $N$ would be normal in $\Lambda$, which is false. Thus $N'N^r$ is not normal in $\Lambda$.

\bigskip

\noindent{\tt G.A.Jones@maths.soton.ac.uk}
\smallskip

\noindent{\tt D.Singerman@maths.soton.ac.uk}
\smallskip

\noindent{\tt paul.watson@psc.ac.uk}

\newpage

\end{document}